\documentclass[11pt]{article}

\usepackage[T2A]{fontenc}
\usepackage[cp1251]{inputenc}
\usepackage{amsthm}%%%%%%%% theorem definition
\usepackage{amsfonts}
\usepackage{eufrak}
\usepackage{amssymb}
\usepackage{amsmath}
\usepackage{cite}


\begin{document}
%%%%%%%%%%%%% begin theorem definition %%%%%%%%%%%%%%%%%%
\newtheoremstyle{mytheorem}
  {\topsep}   % ABOVESPACE
  {\topsep}   % BELOWSPACE
  {\itshape}  % BODYFONT
  {}       % INDENT (empty value is the same as 0pt)
  {\bfseries} % HEADFONT
  {.}         % HEADPUNCT
  {5pt plus 1pt minus 1pt} % HEADSPACE
  { }          % CUSTOM-HEAD-SPEC
\newtheoremstyle{myremark}
  {\topsep}   % ABOVESPACE
  {\topsep}   % BELOWSPACE
  {\upshape}  % BODYFONT
  { }       % INDENT (empty value is the same as 0pt)
  {\bfseries} % HEADFONT
  {.}         % HEADPUNCT
  {5pt plus 1pt minus 1pt} % HEADSPACE
  { }          % CUSTOM-HEAD-SPEC
\theoremstyle{mytheorem}
\newtheorem{theorem}{Theorem}[section]
 \newtheorem{theorema}{Theorem}
 \newtheorem*{a*}{Theorem A}
 \newtheorem*{b*}{Theorem B}
\newtheorem{proposition}{Proposition}[section]
\newtheorem{lemma}{Lemma}[section]
\newtheorem{corollary}{Corollary}[section]
\newtheorem{definition}{Definition}[section]
\theoremstyle{myremark}
\newtheorem{remark}{Remark}[section]
%%%%%%%%%%%%%%%%%%%%% end theorem definition %%%%%%%%%%%%%%%%%%

\centerline{\textbf{The Skitovich--Darmois and Heyde  theorems  for}}
\centerline{\textbf{complex and quaternion random variables}}

\bigskip

\centerline{\textbf{G. M. Feldman}}

\bigskip

 \makebox[20mm]{ }\parbox{125mm}{ \small{We prove the following analogue of the classical Skitovich--Darmois  theorem for complex   random variables. Let $\alpha=a+ib$ be a nonzero complex number.  Then the following statements hold.  $1$. Let either $b\ne 0$, or $b=0$  and $a>0$. Let $\xi_1$ and $\xi_2$ be independent complex random variables.  Assume that the linear forms  $L_1=\xi_1+\xi_2$ and $L_2=\xi_1+\alpha\xi_2$ are independent. Then $\xi_j$ are degenerate random variables.
$2$. Let $b=0$ and $a<0$. Then there exist
complex Gaussian random variables in the wide sense $\xi_1$ and $\xi_2$ such that they are not  complex Gaussian random variables in the narrow sense,  whereas the linear forms $L_1=\xi_1+\xi_2$ and $L_2=\xi_1+\alpha\xi_2$ are independent. We also study an analogue of the Heyde  theorem for complex   random variables.}}

\bigskip


\bigskip


\bigskip

\section{Introduction}

Let $\eta_1$ and $\eta_2$ be real random variables.
Consider the complex random variable $\xi=\eta_1+i\eta_2$.
The random vector $(\eta_1,\eta_2)$ is assigned to it, and the distribution
 $\mu$ of the complex random variable $\xi$ is the distribution of the random vector $(\eta_1,\eta_2)$.
Let $\alpha=a+ib$ be a complex number. Then the random vector $(a\eta_1-b\eta_2, b\eta_1+a\eta_2)$ is assigned to the complex
  random variable  $\alpha\xi$.
Associate with the complex number  $\alpha$ the matrix
\begin{equation}\label{1}
\alpha \longleftrightarrow\left(\begin{matrix}a&-b\\
b&a\end{matrix}\right).
\end{equation}
We denote by $\alpha$ both the complex number $a+ib$ and the corresponding matrix of the form (\ref{1}).
Let $\xi=\eta_1+i\eta_2$ be a complex random variable with distribution $\mu$. Following   \cite{VK1}, we say that $\xi$ is a complex Gaussian random variable
in the wide sense if the distribution of the random vector $(\eta_1,\eta_2)$ is a Gaussian distribution in $ \mathbb{R}^2$.
In this case the characteristic function of the distribution $\mu$ is of the form
\begin{equation}\label{2}
\hat\mu(y)=\exp\{i\langle x, y\rangle-\langle Ay,
y\rangle\},\quad y\in \mathbb{R}^2,
\end{equation}
where $x\in \mathbb{R}^2$, $\langle ., .\rangle$ is the scalar product, and  $A$ is a symmetric positive semidefinite  $(2 \times 2)$-matrix.
If in
$(\ref{2})$  $A$ is a scalar matrix, then we say that the complex random variable
$\xi$ is a complex Gaussian random variable
in the narrow sense.

We note,  on the one hand, that the following result holds. Let $\xi_j$, $j=1, 2, \dots, n$, $n\ge 2$ be independent complex random variables, $\alpha_j$, $\beta_j$ be nonzero complex numbers. Then the independence of the linear forms  $L_1=\alpha_1\xi_1+\dots +\alpha_n\xi_n$
and $L_2=\beta_1\xi_1+\dots +\beta_n\xi_n$ implies that  $\xi_j$ are complex Gaussian random variables in the wide sense. This result is an analogue of the well-known Skitovich--Darmois theorem (\cite[\S\,3.1]{KaLiRa})  for complex Gaussian random variables.
It follows directly from  a weak variant of the Ghurye--Olkin theorem (see Lemma \ref{le2} below), where the Gaussian distribution in the space $\mathbb{R}^m$ is characterized by the independence of two linear forms of $n$ independent random vectors. Coefficients of the linear forms are non-singular matrices. This fact has been noticed, e.g. in   \cite{EK}. On the other hand, the following  theorem, where the complex Gaussian random variables in the narrow sense are characterized by the independence of two linear forms, has been proved in  \cite{VC1}.

\begin{a*} \label{thA} Let $\xi_1$ and $\xi_2$ be independent complex random
variables. Let $\alpha_1, \alpha_2, \beta_1, \beta_2$ be nonzero complex numbers such that
either $\bar \alpha_1 \alpha_2  \beta_1 \bar \beta_2$
or $\sqrt {-\frac{\alpha_1\beta_1}{\alpha_2\beta_2}} \left(\frac{\beta_2}{\beta_1}-\frac{\alpha_2}{\alpha_1}\right)$
is   not a real number. If the linear forms
$L_1=\alpha_1\xi_1+\alpha_2\xi_2$ and  $L_2=\beta_1\xi_1+\beta_2\xi_2$ are
independent, then  $\xi_1$ and $\xi_2$ are complex Gaussian
random
variables in the narrow sense.
\end{a*}

The  proof of Theorem  A, given in \cite{VC1}, is based on
 the Polya characterization theorem for complex Gaussian random variables in
 the narrow sense (see \cite{V2}).
We prove in this note the following statement.

\begin{theorem}\label{th1} Let $\alpha=a+ib$ be a nonzero complex number.  Then the following statements hold.

$1$. Assume that either $b\ne 0$  or $b=0$  and $a>0$. Let $\xi_1$ and $\xi_2$ be independent complex random variables.  Assume that the linear forms  $L_1=\xi_1+\xi_2$ and $L_2=\xi_1+\alpha\xi_2$ are independent. Then $\xi_j$ are degenerate random variables.

$2$. Assume that $b=0$ and $a<0$. Then there exist
complex Gaussian random variables in the wide sense $\xi_1$ and $\xi_2$ such that they are not  complex Gaussian random variables in the narrow sense,  whereas the linear forms $L_1=\xi_1+\xi_2$ and $L_2=\xi_1+\alpha\xi_2$ are independent.
\end{theorem}

Theorem \ref{th1} means that excluding the degenerate case we can not characterize complex Gaussian random variables  in the narrow  sense by the independence of two linear forms of two independent complex random variables.

Note  that on the one hand, a large number of studies has been  devoted to characterization problems of mathematical statistics on different classes of locally compact Abelian groups (see e.g. \cite{SMZh, 2Fe, F1, Fe30, FG,  My1, MiFe1}). On the other hand,  characterization problems of mathematical statistic for complex and quaternion random variables have hardly been studied.

\section{Proof of Theorem \ref{th1}}

To prove Theorem  \ref{th1} we need two lemmas.

\begin{lemma}\label{le1} Let $\xi_j$, $j = 1, 2,\dots, n, \ n \ge 2,$ be independent complex random variables with distributions
 $\mu_j$. Let  $\alpha_j, \beta_j$ be nonzero complex numbers. The linear forms $L_1=\alpha_1\xi_1+\dots+\alpha_n\xi_n$ and $L_2=\beta_1\xi_1+\dots+\beta_n\xi_n$ are independent if and only if the characteristic functions
$\hat\mu_j(y)$ satisfy the equation
\begin{equation}\label{3}
\prod_{j=1}^n\hat\mu_j(\bar \alpha_j u+\bar \beta_j v)
=\prod_{j=1}^n\hat\mu_j(\bar \alpha_j u)\prod_{j=1}^n\hat\mu_j(\bar \beta_j v), \quad u, v \in \mathbb{R}^2,
\end{equation}
where $\bar \alpha_j$ and $\bar \beta_j$ are matrices of the form $(\ref{1})$,
corresponding to the complex numbers $\bar \alpha_j$ and $\bar \beta_j$.
\end{lemma}

The proof of Lemma \ref{le1} is standard. Lemma \ref{le1} is a particular case of a general
statement that concerns to arbitrary locally compact Abelian groups  (see \cite[Lemma 10.1  ]{F1}).

\begin{lemma}\label{le2} Let  $\xi_j$, $j = 1, 2,\dots, n, \ n \ge 2,$ be independent random vectors in the space $\mathbb{R}^m$.
Let
$\beta_j$ be non-singular  $(m\times m)$-matrices
satisfying the conditions
\begin{equation}\label{4}
\beta_i-\beta_j \ \   \mbox{is a non-singular matrix for all} \ i\ne j.
\end{equation}
Then the independence of the linear forms       $L_1 = \xi_1 + \dots
+ \xi_n$ and $L_2 = \beta_1\xi_1 + \dots + \beta_n\xi_n$
implies that all random vectors $\xi_j$  are Gaussian.
\end{lemma}

Lemma \ref{le2} is a weak variant of   the Ghurye--Olkin theorem. In fact  the Ghurye--Olkin theorem states that Lemma \ref{le2} is valid without restriction   (\ref{4}). We note that the main part
of the proof of the Ghurye--Olkin theorem    is the passage from the case when (\ref{4}) is valid to the general case. As to Lemma \ref{le2}, its proof is exactly as
 the proof of the Skitovich--Darmois theorem by the finite difference method (see e.g. \cite[\S 3.2  ]{KaLiRa}).

{\bf Proof of Theorem \ref{th1}}. 1. Denote by $\mu_j$  the distribution of the complex random variable  $\xi_j$. By Lemma \ref{le1}, the characteristic functions
$\hat\mu_j(y)$ satisfy equation (\ref{3}) which takes the form
\begin{equation}\label{5}
\hat\mu_1(u+v)\hat\mu_2(u+\bar \alpha v)
=\hat\mu_1(u)\hat\mu_2(u)\hat\mu_1(v)\hat\mu_2(\bar \alpha v), \quad u, v \in \mathbb{R}^2,
\end{equation}
where $\bar \alpha$ is a matrix of the form $(\ref{1})$,
 corresponding to the complex number $\bar\alpha$. It is obvious that if  $\alpha=1$, then $\mu_j$ are degenerate distributions. So, assume that
    $\alpha\ne 1$, i.e. $1-\alpha\ne 0$. Then it follows from $(\ref{1})$ that condition $(\ref{4})$ holds. Thus, by Lemma \ref{le2}, $\xi_j$ are complex Gaussian random variables in the wide sense. Hence, the characteristic functions $\hat\mu_j(y)$ are of the form
\begin{equation}\label{6}
\hat\mu_1(y)=\exp\{i\langle x_1, y\rangle-\langle Ay,
y\rangle\},\quad \hat\mu_2(y)=\exp\{i\langle x_2, y\rangle-\langle By,
y\rangle\}, \quad y\in \mathbb{R}^2,
\end{equation}
where $x_1, x_2\in \mathbb{R}^2$,
$A=(a_{ij})$, $B=(b_{ij})$ are symmetric positive semidefinite  $(2 \times 2)$-matrices.
Substituting (\ref{6}) into (\ref{5}) we get that the equality
\begin{equation}\label{7}
A+B\bar\alpha=0
\end{equation}
is valid.  Since the matrices $A$ and $B$ are symmetric, it follows from  (\ref{7}) that
$$
b_{11}b+b_{12}a=b_{12}a-b_{22}b.
$$
This implies that
\begin{equation}\label{8}
(b_{11}+b_{22})b=0.
\end{equation}
Since  $B$ is a symmetric positive semidefinite matrix, we have
\begin{equation}\label{9}
b_{11} \ge 0, \quad b_{22} \ge 0, \quad b_{11}b_{22}-b_{12}^2\ge 0.
\end{equation}

Assume that $b\ne 0$. Then it follows from  (\ref{8}) and (\ref{9}) that $b_{11}=b_{22}=b_{12}=0,$ i.e. $B=0$. Then (\ref{7}) implies that $A=0$.
Thus we proved that $\mu_j$ are degenerate distributions.

Assume that $b=0$ and  $a>0$. It follows from (\ref{7}) that the   the equality $A+aB=0$ is valid.
Since $A$ and $B$ are symmetric positive semidefinite matrices, and $a>0$, this implies that $A=B=0$.
Thus, we proved that in this case $\mu_j$ are also degenerate distributions.

2. Assume that $b=0$ and $a<0$. Let $B$ be an arbitrary nonscalar
symmetric positive semidefinite $(2 \times 2)$-matrix. Put $A=-a B$. Let $\xi_1$ and $\xi_2$ be independent complex random variables with distributions $\mu_1$
and $\mu_2$, having the characteristic functions of the form
\begin{equation}\label{10}
\hat\mu_1(y)=\exp\{-\langle Ay,
y\rangle\},\quad \hat\mu_2(y)=\exp\{-\langle By,
y\rangle\}, \quad y\in \mathbb{R}^2.
\end{equation}
Then $\xi_1$ and  $\xi_2$ are complex Gaussian random variables in the wide sense such that they are not  complex Gaussian random variables in the narrow sense.
Since (\ref{7}) is fulfilled, it is easy to that the characteristic functions
$\hat\mu_j(y)$ satisfy equation (\ref{5}). By Lemma \ref{le1},
the linear forms $L_1=\xi_1+\xi_2$ and  $L_2=\xi_1+\alpha\xi_2$ are independent. $\Box$

\begin{remark}\label{r1} Let $\xi_1$ and $\xi_2$ be independent complex random
variables with distributions $\mu_1$ and $\mu_2$. Let $\alpha_1, \alpha_2, \beta_1, \beta_2$ be nonzero complex numbers. Consider the linear forms $L_1=\alpha_1\xi_1+\alpha_2\xi_2$ and  $L_2=\beta_1\xi_1+\beta_2\xi_2$ and assume that $L_1$ and $L_2$ are independent. We want to describe the possible distributions  $\mu_1$
and $\mu_2$. Introduce into consideration new independent complex random variables $\xi_1'=\alpha_1\xi_1$ and $\xi_2'=\alpha_2\xi_2$ and note that the linear forms $L_1$ and  $L_2$ are independent if and only if the linear forms $L_1$ and  $cL_2$ are independent  for any nonzero complex $c$. From this it follows that the description of possible distributions  $\mu_1$
and $\mu_2$ is reduced to the case when the linear forms $L_1$ and  $L_2$ are of the form
$L_1=\xi_1+\xi_2$ and $L_2=\xi_1+\alpha\xi_2$, where $\alpha=\alpha_1\alpha_2^{-1}\beta_1^{-1}\beta_2$, i.e. to Theorem \ref{th1}.
\end{remark}
\begin{remark}\label{r2} As has been noted above Theorem \ref{th1} implies that excluding the degenerate case we can not characterize complex Gaussian random variables  in the narrow  sense by the independence of two linear forms $L_1$ and  $L_2$ of two independent complex random variables $\xi_1$ and $\xi_2$.
The situation will not change if we   consider
 $n$ independent complex random variables $\xi_j$. The following proposition shows that do not exist   coefficients $\alpha_j, \beta_j$ such that the independence of the linear forms  $L_1=\alpha_1\xi_1+\dots+\alpha_n\xi_n$ and  $L_2=\beta_1\xi_1+\dots+\beta_n\xi_n$ implies that $\xi_j$ are complex Gaussian random variables in the narrow sense, but need not be degenerate.
\end{remark}
\begin{proposition}\label{pr1}  Let $\xi_j$, $j = 1, 2,\dots, n, \ n \ge 2,$ be independent nondegenerate complex Gaussian random variables in the narrow sense. Let
$\beta_j$ be nonzero complex numbers. Assume that the linear forms   $L_1=\xi_1+\dots+\xi_n$ and  $L_2=\beta_1\xi_1+\dots+\beta_n\xi_n$ are independent.
Then there exist independent   complex Gaussian random variables in the wide sense $\eta_j$ such that they are not  complex Gaussian random variables in the narrow sense, whereas the linear forms  $L'_1=\eta_1+\dots+\eta_n$ and  $L'_2=\beta_1\eta_1+\dots+\beta_n\eta_n$ are also independent.
\end{proposition}

{\bf Proof}. Denote by $\mu_j$ the distribution of the complex Gaussian random variables $\xi_j$. By the condition, the characteristic functions $\hat\mu_j(y)$ are of the form
\begin{equation}\label{11}
\hat\mu_j(y)=\exp\{i\langle x_j, y\rangle-\sigma_j\langle y,
y\rangle\},\quad y\in \mathbb{R}^2,
\end{equation}
where $x_j\in \mathbb{R}^2$, $\sigma_j>0$. By Lemma \ref{le1}, the characteristic functions
$\hat\mu_j(y)$ satisfy equation  (\ref{3}). Substituting (\ref{11})
into (\ref{3}), we find
\begin{equation}\label{12}
\sum_{j=1}^n\sigma_j\bar \beta_j=0,
\end{equation}
where $\bar \beta_j$  are matrices of the form $(\ref{1})$
corresponding to the complex numbers $\bar \beta_j$. Let $A$ be an arbitrary
nonscalar  symmetric positive semidefinite  $(2 \times 2)$-matrix. Put $A_j=\sigma_j A$. Let $\eta_j$ be independent complex random variables with distributions  $\nu_j$ such that their characteristic functions are of the form
\begin{equation}\label{13}
\hat\nu_j(y)=\exp\{-\langle A_jy,
y\rangle\}, \quad y\in \mathbb{R}^2.
\end{equation}
Then $\eta_j$ are complex Gaussian random variables in the wide sense such that they are not  complex Gaussian random variables in the narrow sense.
 Substitute (\ref{13}) into (\ref{3}). It is clear that the characteristic functions   $\hat\nu_j(y)$ satisfy equation (\ref{3}) if and only if the equality
\begin{equation}\label{14}
\sum_{j=1}^nA_j\bar \beta_j=0
\end{equation}
holds. It is obvious that (\ref{14}) follows from (\ref{12}). By Lemma \ref{le1},
the linear forms $L'_1=\eta_1+\dots+\eta_n$ and  $L'_2=\beta_1\eta_1+\dots+\beta_n\eta_n$ are independent. $\Box$

\section{The Heyde theorem  for the complex random variables}

The statement closely connected with the Skitovich--Darmois theorem was proved by Heyde (\cite{He}, see
also \cite[\S\,13.4.1]{KaLiRa}). According to Heyde's theorem the Gaussian distribution
on the real line is characterized by the symmetry of the conditional
distribution of one linear form of
$n$ independent random variables given another. For $n=2$,
  this theorem can be formulated as follows.

\begin{b*} Let $\xi_1$  and $\xi_2$ be independent random variables.
 Let $a$ be a nonzero real number, $a\ne -1$.
If the conditional distribution of  the linear form  $L_2 = \xi_1 +
a\xi_2$  given $L_1 = \xi_1 + \xi_2$ is symmetric, then
the random variables  $\xi_j$   are Gaussian.
\end{b*}

Theorem \ref{th1} allows us to prove an analogue of Theorem B for complex random variables.

\begin{theorem}\label{th2} Let $\alpha=a+ib$ be a nonzero complex number,
  $\alpha\ne -1$. Then the following statements hold.

$1$. Assume that either $b\ne 0$  or $b=0$  and $a>0$. Let $\xi_1$ and $\xi_2$ be independent complex random variables.  Assume that the conditional distribution of  the linear form  $L_2 = \xi_1 +
\alpha\xi_2$  given $L_1 = \xi_1 + \xi_2$ is symmetric. Then $\xi_j$ are degenerate random variables.

$2$. Assume that  $b= 0$ and $a<0$.
Then there exist
complex Gaussian random variables in the wide sense $\xi_1$ and $\xi_2$ such that they are not  complex Gaussian random variables in the narrow sense,  whereas the conditional distribution of  the linear form  $L_2 = \xi_1 +
\alpha\xi_2$  given $L_1 = \xi_1 + \xi_2$ is symmetric.
\end{theorem}

To prove Theorem \ref{th2} we need two lemmas. The following lemma was proved in  \cite{My1}
for random variables with values in a locally compact Abelian group.
We formulate it for complex random variables.

\begin{lemma}\label{le3}  Let  $\xi_1$ and  $\xi_2$ be independent complex random variables.
Let $\alpha$ be a nonzero complex number.   If the conditional distribution of the linear form    $L_2 = \xi_1 +
\alpha\xi_2$  given $L_1 = \xi_1 +
\xi_2$  is symmetric,   then the linear forms
$M_1=(I+\alpha)\xi_1+2\alpha\xi_2$ and
$M_2=2\xi_1+(I+\alpha)\xi_2$  are independent.
\end{lemma}

\begin{lemma}\label{le4} Let  $\xi_1$ and  $\xi_2$ be independent complex random variables
with distributions $\mu_1$ and $\mu_2$.
 Let $\alpha$ be a nonzero complex number.  The conditional distribution of the linear form
 $L_2 = \xi_1 + \alpha\xi_2$
 given $L_1 = \xi_1 + \xi_2$ is symmetric if and only
 if the characteristic functions
 $\hat\mu_j(y)$ satisfy the equation
\begin{equation}
\label{15} \hat\mu_1( u+ v )\hat\mu_2( u+\bar\alpha
v )= \hat\mu_1( u- v )\hat\mu_2( u-\bar\alpha v ), \quad u, v\in \mathbb{R}^2,
\end{equation}
where $\bar \alpha$ is a matrix of the form $(\ref{1})$
corresponding  to the complex number $\bar\alpha$.
\end{lemma}

Lemma \ref{le4} is a particular case of a general
statement that concerns to locally compact Abelian groups
(see \cite[Lemma 16.1  ]{F1}).

{\bf Proof of Theorem \ref{th2}}. 1. By Lemma \ref{le3}, the linear forms $M_1=(1+\alpha)\xi_1+2\alpha\xi_2$ and
$M_2=2\xi_1+(1+\alpha)\xi_2$ are independent. Putting $\xi_1'=(1+\alpha)\xi_1$,
 $\xi_2'=2\alpha\xi_2$, we obtain that the linear forms $N_1=\xi_1'+\xi_2'$ and
$N_2=\frac{2}{1+\alpha}\xi_1'+\frac{1+\alpha}{2\alpha}\xi_2'$ are independent.
Hence, the linear forms
 $P_1=\xi_1'+\xi_2'$ and
$P_2=\xi_1'+\beta\xi_2'$, where $\beta=\frac{(1+\alpha)^2}{4\alpha}$, are also independent.
We have
\begin{equation}
\label{16n} \textstyle \beta=\frac{1}{4}\left(\left(a+2+\frac{a}{|\alpha|^2}\right)+ib\left(1-
\frac{1}{|\alpha|^2}\right)\right)=p+iq.
\end{equation}
Taking into account $(\ref{16n})$ it is easy to verify that the following statements are valid.

A. If $b\ne 0$ and   $|\alpha|\ne 1$, then  $q\ne 0$.

B. If $|\alpha|= 1$, then $q=0$ and $p>0$.

C. If $b=0$, then $q=0$ and  $p$ and $a$ have the same signs.

Assume that either $b\ne 0$, or $b=0$  and $a>0$. Taking into account statements A--C, apply Theorem \ref{th1} to the independent
random variables $\xi_1'$ and $\xi_2'$ and to the linear forms $P_1$ and $P_2$.  We obtain
 that $\xi_j'$ are degenerate random variables. Then, obviously,  $\xi_j$
are also degenerate random variables.

2. Assume that $b=0$ and $a<0$. It is easy to verify that condition
 $(\ref{7})$ is necessary and sufficient in order that characteristic functions
 of the form $(\ref{10})$ satisfy equation  $(\ref{15})$. Reasoning
 as in the proof of case 2 in Theorem \ref{th1}, and applying Lemma \ref{le4}, instead of Lemma \ref{le1}, we complete the proof of Theorem \ref{th2}. $\Box$

Note that
Theorem \ref{th2} implies that  excluding the degenerate case we can not characterize complex Gaussian random variables  in the narrow  sense by the symmetry of the conditional distribution of  one linear form  of two independent random variables  given another.

\begin{remark}\label{r3} Let $\xi_1$ and $\xi_2$ be independent complex random
variables with distributions $\mu_1$ and $\mu_2$. Let $\alpha_1, \alpha_2, \beta_1, \beta_2$ be nonzero complex numbers. Consider the linear forms $L_1=\alpha_1\xi_1+\alpha_2\xi_2$ and  $L_2=\beta_1\xi_1+\beta_2\xi_2$ and assume that the conditional distribution of  the linear form $L_2$ given  $L_1$ is symmetric. It is easy to see that the description of possible distributions  $\mu_1$
and $\mu_2$ is reduced to the case when the linear forms $L_1$ and  $L_2$ are of the form
$L_1=\xi_1+\xi_2$ and $L_2=\xi_1+\alpha\xi_2$, where $\alpha=\alpha_1\alpha_2^{-1}\beta_1^{-1}\beta_2$, i.e. to Theorem \ref{th2}.
\end{remark}

\section{Characterization theorems
 for  quaternion random variables}

Let $\eta_1$, $\eta_2$, $\eta_3$, $\eta_4$ be real random variables.
Consider the quaternion  random variable
 $\xi=\eta_1+i\eta_2+j\eta_3+k\eta_4$.
The random vector $(\eta_1,\eta_2,\eta_3,\eta_4)$ is assigned to it, and the distribution
 $\mu$ of the quaternion random variable $\xi$ is the distribution of the random vector $(\eta_1,\eta_2,\eta_3,\eta_4)$.
Let $\alpha=a+ib+jc+kd$ be a quaternion. Then the random vector $(a\eta_1-b\eta_2-c\eta_3-d\eta_4, b\eta_1+a\eta_2-d\eta_3+c\eta_4, c\eta_1+d\eta_2+a\eta_3-b\eta_4, d\eta_1-c\eta_2+b\eta_3+a\eta_4)$ is assigned to  the quaternion
  random variable  $\alpha\xi$.
Associate with the quaternion  $\alpha$ the matrix
\begin{equation}\label{16}
\alpha \longleftrightarrow\left(\begin{matrix}a&-b&-c&-d\\
b&a&-d&c\\
c&d&a&-b\\
d&-c&b&a\end{matrix}\right).
\end{equation}
We denote by $\alpha$ both the quaternion $a+ib+jc+kd$ and the corresponding matrix of the form (\ref{16}). The quaternion  Gaussian random variables
in the wide sense and in the narrow sense are defined in the same way as  for complex Gaussian random variables (see \cite{V, VC2}).
A theorem similar to Theorem A was proved in   \cite{VC1} for quaternion  random variables. We shall prove the following statement.

\begin{theorem}\label{th3} Let $\alpha=a+ib+jc+kd$ be a nonzero quaternion. Then the following statements hold.

$1$. Assume that either the imaginary part of $\alpha$ is not equal to zero  or the imaginary part of $\alpha$ is equal to zero  and $a>0$. Let $\xi_1$ and $\xi_2$ be independent quaternion random variables.  Assume that the linear forms  $L_1=\xi_1+\xi_2$ and $L_2=\xi_1+\alpha\xi_2$ are independent. Then $\xi_j$ are degenerate random variables.

$2$. Assume that the imaginary part of $\alpha$ is equal to zero  and $a>0$. Then there exist
quaternion Gaussian random variables in the wide sense $\xi_1$ and $\xi_2$ such that they are not  quaternion Gaussian random variables in the narrow sense,  whereas the linear forms $L_1=\xi_1+\xi_2$ and $L_2=\xi_1+\alpha\xi_2$ are independent.
\end{theorem}

{\bf Proof}. We follow the scheme of the proof of Theorem \ref{th1}. Denote by $\mu_j$  the distribution of the quaternion  random variable  $\xi_j$. A lemma similar to Lemma \ref{le1} is valid for quaternion random variables. It implies that the characteristic functions
$\hat\mu_j(y)$ satisfy the equation
\begin{equation}\label{17}
\hat\mu_1(u+v)\hat\mu_2(u+\bar \alpha v)
=\hat\mu_1(u)\hat\mu_2(u)\hat\mu_1(v)\hat\mu_2(\bar \alpha v), \quad u, v \in \mathbb{R}^4,
\end{equation}
where $\bar \alpha$ is a matrix of the form $(\ref{16})$,
corresponding to the quaternion $\bar\alpha$. It is obvious that if  $\alpha=1$, then $\mu_j$ are degenerate distributions. So, assume that
    $\alpha\ne 1$, i.e. $1-\alpha\ne 0$. Then it follows from $(\ref{16})$ that condition $(\ref{4})$ holds. Thus, by Lemma \ref{le2}, $\xi_j$ are quaternion Gaussian random variables in the wide sense. Hence, the characteristic functions
$\hat\mu_j(y)$ are of the form
\begin{equation}\label{new1}
\hat\mu_1(y)=\exp\{i\langle x_1, y\rangle-\langle Ay,
y\rangle\},\quad \hat\mu_2(y)=\exp\{i\langle x_2, y\rangle-\langle By,
y\rangle\}, \quad y\in \mathbb{R}^4,
\end{equation}
where $x_1, x_2\in \mathbb{R}^4$,
$A=(a_{ij})$, $B=(b_{ij})$ are symmetric positive semidefinite $(4 \times 4)$-matrices.
Substituting  (\ref{new1}) into (\ref{17}), we get that equality  (\ref{7}) holds. Since the matrices $A$ and $B$ are symmetric,  (\ref{7}) implies that  the numbers $b_{ij}, a, b, c$ and $d$ satisfy a system of equations, and we find from it
\begin{equation}\label{18}
\begin{cases}\ b(b_{11}+b_{22}+b_{33}+b_{44})=0\\ \ c(b_{11}+b_{22}+b_{33}+b_{44})=0\\
\ d(b_{11}+b_{22}+b_{33}+b_{44})=0\\
\end{cases}
\end{equation}
It follows from (\ref{18}) that if the imaginary part of the quaternion $\alpha$ is not equal to zero, then
$b_{11}=b_{22}=b_{33}=b_{44}=0$. Since $B$ is a symmetric positive semidefinite
matrix, this implies that $B=0,$ and in view of (\ref{7}),  $A=0.$ It means that $\xi_j$ are degenerate random variables.
In the case if the imaginary part of the quaternion $\alpha$ is   equal to zero, we reason exactly as in the consideration of the corresponding case in Theorem \ref{th1}.  $\Box$

We note that Remarks \ref{r1} and \ref{r2}, and Proposition \ref{pr1} are also valid for quaternion random variables. The proof of Proposition \ref{pr1} in the quaternion case is the same as in the complex one.

Theorem \ref{th3} allows us also to prove an analogue of Theorem B for quaternion random variables.

\begin{theorem}\label{th4} Let $\alpha=a+ib+jc+kd$ be a nonzero quaternion,
  $\alpha\ne -1$. Then the following statements hold.

  $1$. Assume that either the imaginary part of $\alpha$ is not equal to zero  or the imaginary part of $\alpha$ is equal to zero  and $a>0$. Let $\xi_1$ and $\xi_2$ be independent quaternion random variables.  Assume that the conditional distribution of  the linear form  $L_2 = \xi_1 +
\alpha\xi_2$  given $L_1 = \xi_1 + \xi_2$ is symmetric. Then $\xi_j$ are degenerate random variables.

$2$. Assume that the imaginary part of $\alpha$ is equal to zero  and  $a<0$.
Then there exist
quaternion Gaussian random variables in the wide sense $\xi_1$ and $\xi_2$ such that they are not  quaternion Gaussian random variables in the narrow sense,  whereas the conditional distribution of  the linear form  $L_2 = \xi_1 +
\alpha\xi_2$  given $L_1 = \xi_1 + \xi_2$ is symmetric.
\end{theorem}

{\bf Proof}. We reason as in the proof of Theorem \ref{th2} and retain the same notation. In so doing, the quaternion  $\beta=\frac{(1+\alpha)^2}{4\alpha}$   is of the form
$$
\textstyle\beta=\frac{1}{4}\left(\left(a+2+\frac{a}{|\alpha|^2}\right)
+(ib+jc+kd)\left(1-\frac{1}{|\alpha|^2}\right)\right),
$$
Lemmas \ref{le3} and \ref{le4} are also valid for quaternion random variables,
and we use Theorem \ref{th3} instead of Theorem \ref{th1}.  $\Box$

\newpage

\noindent Gennadiy Feldman

\medskip

\noindent B. Verkin Institute for Low Temperature Physics and Engineering \\ of the National Academy
of Sciences of Ukraine\\ Nauky Ave. 47,\\
Kharkiv, 103,\\ 61103, Ukraine

\medskip

\noindent e-mail: feldman@ilt.kharkov.ua

\end{document}